\documentclass[11pt]{article}
\usepackage{amsmath,amssymb}
\newcommand{\doublespace}
   {\addtolength{\baselineskip}{0.15\baselineskip}}

\newcommand{\qed}{\hfill\rule{2mm}{2mm}\\}
\newtheorem{pdef}{Definition}[section] %
\newtheorem{thm}[pdef]{Theorem}        
\newtheorem{cor}[pdef]{Corollary}
\newtheorem{lem}[pdef]{Lemma}

\newcounter{equationnumber}
\renewcommand{\theequation}{\thesection.\arabic{equation}}
\def\mathletters{
    \addtocounter{equation}{1}
    \edef\@currentlabel{\theequation}
    \setcounter{equationnumber}{\value{equation}}
    \setcounter{equation}{0}
    \edef\theequation{\@currentlabel\noexpand\alph{equation}}
    }

\title{Characterization of the matrix whose norm is determined by
its action on decreasing sequences}
\author{Chang-Pao Chen$^\ast$, Hao-Wei Huang, and Chun-Yen Shen\\
\vspace{0.1in}\\
Department of Mathematics\\
National Tsing Hua University\\
Hsinchu, Taiwan 300\\
R. O. C.\\
Email: cpchen@math.nthu.edu.tw}
\date{ }
\begin{document}
\maketitle
\doublespace
\pagestyle{myheadings} \thispagestyle{plain} \markboth{   }{ }
\begin{abstract}
Let $A=(a_{j,k})_{j,k \ge 1}$ be a non-negative matrix. In this
paper, we characterize those $A$ for which $\|A\|_{E, F}$ are
determined by their actions on decreasing sequences, where $E$ and
$F$ are suitable normed Riesz spaces of sequences. In particular,
our results can apply to the following spaces: $\ell_p$, $d(w,p)$,
and $\ell_p(w)$. The results established here generalize the
corresponding ones given by Bennett in Quart. J. Math. Oxford (2),
49(1998), 395-432, by Chen et al in J. Math. Anal. Appl.
273(2002), 160-171 and by Jameson in Illinois J. Math. 43(1999),
79-99.
\end{abstract}
\footnotetext[1]{This work is supported by the National Science
Council, Taipei, ROC, under Grant NSC 94-2115-M-007-008}
\footnotetext[2]{{\it 2000 Mathematics Subject Classification:}\,
Primary 15A60, 40G05, 47A30, 47B37, Secondary 46B42.}
\footnotetext[3]{{\it Key words and phrases.}\, norms of matrices,
normed Riesz spaces, weighted mean matrices, N\"orlund mean
matrices, summability matrices, matrices with row decreasing.}

\section{Introduction}

Let $w_1\ge w_2\ge\cdots\ge 0$. For $1\le p\le\infty$, denote by
$\ell_p(w)$ the space of all sequences $x=\{x_k\}_{k=1}^\infty$
such that
$$\displaystyle
 \|x\|_{\ell_p(w)}:=\left(\sum_{k=1}^\infty
|x_k|^pw_k\right)^{1/p}<\infty.\leqno{(1.1)}$$ The Lorentz
sequence space $d(w,p)$ is the space of null sequences $x$ for
which $x^*$ is in $\ell_p(w)$, with norm
$\|x\|_{w,p}=\|x^*\|_{\ell_p(w)}$, (cf. [1, 7]).  Here $x^*$ is
the decreasing rearrangement of $\{|x_k|\}_{k=1}^\infty$. When
$w_k=1$ for all $k$, $\ell_p(w)$ coincides with $\ell_p$ in the
usual sense (the norm of which we denote by $\|\cdot\|_p$). We
also have $\ell_\infty(w)=\ell_\infty$ for any $w$. We write $x\ge
0$ if $x_k\ge 0$ for all $k$. Similarly, $x\downarrow$ will mean
that $\{x_k\}_{k=1}^\infty$ is decreasing, that is, $x_k\ge
x_{k+1}$ for all $k\ge 1$. For a non-negative matrix
$A=(a_{j,k})_{j,k\ge 1}$ and two normed sequence spaces $(E,
\|\cdot\|_E)$, $(F, \|\cdot\|_F)$ in $\ell_p(w)$,  let
$\|A\|_{E,F}$ denote the norm of $A$ when regarded as an operator
from $E$ to $F$. Clearly, for $A \geq 0$, the norm of A is
determined by non-negative sequeances and $\|A\|_{E,F}\ge
\|A\|_{E,F,\downarrow}$, where $$ \|A\|_{E,
F,\downarrow}:=\sup_{\|x\|_E=1,x\ge 0, x\downarrow}\|Ax\|_F. $$ In
[3, Problem 7.23], Bennett asked the following question for
$E=F=\ell_p$: When does the equality
$\|A\|_{E,F}=\|A\|_{E,F,\downarrow}$ hold ? It is one of great
importance in the general theory of inequalities.

In [2, page 422] and [3, page 422], Bennett  established this
upper bound equality for the case that $E=F=\ell_p$, $1<p<\infty$,
and $A$ is a weighted mean matrix with decreasing weights $w_n$.
This result was extended by Jameson [6, Theorem 2] to the case
that $E=F$ is a Banach lattice of sequences with property $(PS)$
and  $A$ satisfies the following condition:
$$
\sum_{j=1}^l\sum_{k=1}^r a_{j,k} \geq\sum_{j\in N_l}\sum_{k\in
N_r}a_{j,k}\qquad (l,r\ge 1; |N_l|=l, |N_r|=r). \leqno{(1.2)}$$
For the definition of $(PS)$, we refer the readers to \S3. Here
$N_s$ denotes a set of positive integers having $s$ elements and
$|N_s|=s$ stands for all possibilities of $N_s$. Later, in a joint
paper, the first present author extended Bennett's result in a
different direction. More precisely, in [5, Lemma 2.4], Chen et al
established the equality $\|A\|_{E,F}=\|A\|_{E,F,\downarrow}$ for
the case that $E=F=\ell_p$, $1<p<\infty$, and  $A$ is a
non-negative lower triangular matrix with rows decreasing in the
sense that $a_{j,k}\ge a_{j,k+1}$ for all $j,k\ge 1$.

The purpose of this paper is to extend the results of Bennett,
Jameson and Chen-Luor-Ou to a more general setting. In \S2, we
introduce the collection ${\cal R}^{\gamma,\lambda}_A$, which is a
special set of row rearrangements of $A$ with indices
$\gamma\le\lambda$. We prove that for a non-negative $n\times
\infty$ matrix $A$, ${\cal R}^{\gamma,\lambda}_A\neq \phi$ for
some pair $(\gamma, \lambda)$ with $0\le \gamma\le\lambda\le n$.
We also prove that for $x=\{x_k\}_{k=1}^\infty\ge 0$, there exists
some $B\in\cup_{0\le\gamma\le\lambda\le n}
 {\cal R}_A^{\gamma,\lambda}$, depending on $A$ and $x$,
such that the finite sequence $Bx=\{\sum_{k=1}^\infty
b_{j,k}x_k\}_{j=1}^n$ is decreasing. Based on these, we establish
in Theorem 3.2 the upper bound equality for the case that $E$ and
$F$ are two suitable normed Riesz spaces of sequences with
property $(PS)$ and the following condition is satisfied by some
positive integer $n_0$:
\begin{itemize}
\item[$(1.3)$] \quad for given $n\ge n_0$ and
$\displaystyle B=(b_{j,k})\in\cup_{0\le\gamma\le\lambda\le n}
  {\cal R}_{A_n}^{\gamma,\lambda}$,
there exists some $C\in {\cal  R}_{A_n}$, depending on $n$ and
$B$, such that the following inequality holds:
$$
\sum_{j=1}^l\sum_{k=1}^r c_{j,k} \ge\sum_{j=1}^l\sum_{k\in
N_r}b_{j,k}\qquad (1\le l\le n;\, r\ge 1;\, |N_r|=r),$$
\end{itemize}
 where $C=(c_{j,k})$, $A_n$ is the $n\times
\infty$ matrix obtained from the first $n$ rows of $A$, and $\cal
R_A$ is the set of all row rearrangements of $A$. In particular,
Theorem 3.2 can apply to any of $\ell_p$ and $d(w, p)$ for the
spaces $E$ and $F$. However, $\ell_\infty$ is excluded. A similar
result is also established for the case $F=\ell_p(w)$, (cf.
Theorem 3.3). In \S4, we shall give a detailed investigation of
$(1.3)$ for the matrix $A$. These include the investigations of
the Hilbert matrix, the weighted mean matrix, the N\"orlund
matrix, summability matrices, and matrices with row decreasing. Of
course, the Gamma matrix $\Gamma(\alpha)$ and the Ces\`aro matrix
$C(\alpha)$ are also examined. Since $(1.2)\Longrightarrow (1.3)$
(by choosing $ C=A_n$), our results generalize [6, Theorem 2] and
Bennett's result. On the other hand, $(1.3)$ is satisfied,
provided $A$ is row decreasing. In this case, we choose $C=B$.
Therefore, our results (especially
 Corollary 4.7) also include [5, Lemma 2.4] as a special case. We
refer the readers to \S4 for details.

\section{The collection ${\cal R}^{\gamma,\lambda}_A$}

Let $A=(a_{j,k})$ be an $n\times \infty$ matrix. Here $1\le j\le
n$ and $1\le k<\infty$. We say that an $n\times \infty$ matrix
$B=(b_{j,k})$ is a matrix obtained from $A$ by row rearrangements,
if there is a one-to-one mapping $\sigma$ from $\{1,2,\cdots,n\}$
onto itself with $b_{j,k}=a_{\sigma(j),k}$ for all $j$ and for all
$k$. Denote  by ${\cal R}_A$ the collection of these matrices.
Clearly, $A\in {\cal R}_A$. We pay attention to the following
subset of ${\cal R}_A$.

\begin{pdef}
For $0\le\gamma\le\lambda\le n$, we write $B\in{\cal
R}_A^{\gamma,\lambda}$ if and only if $B\in{\cal R}_A$ and
$B=(b_{j,k})$ is of the form:
\begin{itemize}
\item[$(i)$] $b_{j,k}\ge b_{j+1,k}$ for $j\le\gamma$ or $j\ge\lambda$,
\item[$(ii)$]   $b_{r_1,k}\ge b_{j,k}\ge b_{r_2,k}$ for
     $r_1\le\gamma<j<\lambda\le r_2,$
\item[$(iii)$] No $\alpha$ with $\gamma<\alpha<\lambda$ possesses the property:
$b_{\alpha,k}\ge b_{j,k}$ for all $\gamma<j<\lambda$ or
$b_{\alpha,k}\le b_{j,k}$ for all $\gamma<j<\lambda$.
\end{itemize}
\end{pdef}

By definition, no row of the matrices in ${\cal R}_A^{0,\lambda}$
 is greater than or equal to the other rows. Analogously, no row of the matrices in ${\cal R}_A^{\gamma,n}$
 is less than
or equal to the other rows. Moreover, each matrix $B=(b_{j,k})$ in
${\cal R}_A^{\lambda,\lambda}$ or ${\cal R}_A^{\lambda,\lambda+1}$
must be column decreasing, that is, $
 b_{1,k}\ge b_{2,k}\ge\cdots
 \ge b_{n,k}
$ for all $k$. For $A$ with column decreasing, $\displaystyle
 B\in \cup_{0\le\gamma\le\lambda\le n}
 {\cal R}_A^{\gamma,\lambda}$ if and only if $B=A$.

\begin{lem}
${\cal R}_A^{\gamma,\lambda}\neq \emptyset$ for some pair
$(\gamma,\lambda)$ with $0\le \gamma\le\lambda\le n$.
\end{lem}

\begin{pf} We shall prove the existence of a matrix
$B=(b_{j,k})$ with $B\in{\cal R}_A^{\gamma,\lambda}$  for some
pair $(\gamma,\lambda)$ obeying the condition $0\le
\gamma\le\lambda\le n$.  Fix a  row $(a_{j_1,1},a_{j_1,2},\cdots)$
of $A$ and check whether $a_{j_1,k}\ge  a_{j,k}$ for all $j$ and
for all $k$ with $j\neq j_1$. We can consider $j_1$ in the order:
$j_1=1,2,\cdots, n$. If so, let
$(b_{1,1},b_{1,2},\cdots)=(a_{j_1,1},a_{j_1,2},\cdots)$ and choose
another row, say $(a_{j_2,1},a_{j_2,2},\cdots)$, from the other
$n-1$ rows. Check whether $a_{j_2,k}\ge a_{j,k}$ for all $j$ and
for all $k$ with $j\neq j_1,j_2$. If so, let
$(b_{2,1},b_{2,2},\cdots)=(a_{j_2,1},a_{j_2,2},\cdots)$ and choose
another row, say $(a_{j_3,1},a_{j_3,2},\cdots)$, from the other
$n-2$ rows. Check whether $a_{j_3,k}\ge a_{j,k}$ for all $j$ and
for all $k$ with  $j\neq j_1,j_2,j_3$. Continue this process up to
the maximal possibility. We shall stop at some step, say the
$\gamma$th step, and we shall find the first $\gamma$ rows of $B$
with the property: $b_{j,k}\geq b_{j+1,k}$ for all $1\le j<\gamma$
and $b_{\gamma,k}\geq a_{j,k}$ for all $j\neq
j_1,j_2,\cdots,j_\gamma$. Apply the same procedure to the
remainder of rows in the following way. First, choose a row, say
$(a_{s_1,1},a_{s_1,2},\cdots)$, and check whether $a_{s_1,k}\leq
a_{j,k}$ for all $j$ and for all $k$ with $j\neq
j_1,j_2\cdots,j_\gamma,s_1$. If so, let
$(b_{n,1},b_{n,2},\cdots)=(a_{s_1,1},a_{s_1,2},\cdots)$ and choose
a new row, say, $(a_{s_2,1},a_{s_2,2},\cdots)$, from the other
$n-\gamma-1$ rows of $A$. Check whether $a_{s_2,k}\leq a_{j,k}$
for all $j$ and for all $k$ with $j\neq
j_1,\cdots,j_\gamma,s_1,s_2$. If so, let
$(b_{n-1,1},b_{n-1,2},\cdots)=(a_{s_2,1},a_{s_2,2},\cdots)$.
Continue this process up to the maximal possibility. We will stop
at  some step, which corresponds to the $\lambda$th row of $B$. We
also find the last $(n-\lambda+1)$ rows of $B$ with the property:
$b_{j,k}\geq b_{j+1,k}$ for all $\lambda\le j<n$  and
$b_{\lambda,k}\le b_{j,k}$ for all $j\neq j_1,j_2,\cdots,j_\gamma,
s_1,s_2,\cdots,s_{n-\lambda+1}.$ Put the rest of rows into the
middle block of $B$ in any order. Then the final matrix $B$ has
the prescribed property. This completes the proof.\qed
\end{pf}

\begin{lem}
Let $A=(a_{j,k})$ be a non-negative $n\times \infty$ matrix and
$x=\{x_k\}_{k=1}^\infty\ge 0$. Then there exists some
$B\in\cup_{0\le\gamma\le\lambda\le n}
 {\cal R}_A^{\gamma,\lambda}$, depending on $A$ and $x$,
such that the sequence $Bx=\{\sum_{k=1}^\infty
b_{j,k}x_k\}_{j=1}^n$ is decreasing.
\end{lem}

\begin{pf} Lemma 2.2 guarantees the existence of a matrix
 $B\in\cup_{0\le\gamma\le\lambda\le n}
 {\cal R}_A^{\gamma,\lambda}$, say $B\in {\cal R}_A^{\gamma,\lambda}$. Let
$y_j=\sum_{k=1}^\infty b_{j,k}x_k$. By Definition 2.1(i), we
obtain $y_1\ge y_2\ge\cdots\ge y_\gamma$ and $y_\lambda\ge
y_{\lambda+1}\ge\cdots\ge y_n$.  From Definition 2.1(ii), we see
that $y_\gamma\ge y_j\ge y_\lambda$ for all $\gamma<j<\lambda.$
Make a decreasing rearrangement for
$y_{\gamma+1},\cdots,y_{\lambda-1}$, and let $\tilde B$ be the
corresponding matrix by applying such row rearrangements to $B$ .
It is clear that $\tilde B$ still lies in the set ${\cal
R}_A^{\gamma,\lambda}$ and it has the prescribed property. We
complete the proof.\qed
\end{pf}

\section{Main Results}

 We have the following result.

\begin{lem}
Let $\{v_k\}_{k=1}^n$ and $\{u_k\}_{k=1}^n$ be two non-negative
sequences such that
$$
 \sum_{k=1}^r v_k\ge\sum_{k\in N_r}u_k
 \qquad(r=1,\ldots,n;\,|N_r|=r).\leqno{(3.1)}
$$
Then
$$
 \sum_{k=1}^nv_kx_k^*\ge\sum_{k=1}^n u_kx_k\qquad(x=\{x_k\}_{k=1}^n\ge 0).
$$
\end{lem}

\begin{pf}
We have $x_k^*-x_{k+1}^*\ge 0$ for all $1\le k<n$. Let $\{\tilde
u_k\}_{k=1}^n$ denote the corresponding rearrangement of
$\{u_k\}_{k=1}^n$ such that $\sum_{k=1}^n u_kx_k=\sum_{k=1}^n
\tilde u_kx_k^*$. Employing the summation by parts and $(3.1)$, we
get
\begin{align*}
 \sum_{k=1}^n u_kx_k&=\sum_{k=1}^n\tilde u_kx_k^*
=\sum_{k=1}^{n-1}(x_k^*-x_{k+1}^*)\biggl (\sum_{s=1}^k\tilde
u_s\biggr )
 +x_n^*\biggl (\sum_{k=1}^n\tilde u_k\biggr )\\
&\le\sum_{k=1}^{n-1}(x_k^*-x_{k+1}^*)\biggl (\sum_{s=1}^kv_s\biggr
)
 +x_n^*\biggl (\sum_{k=1}^nv_k\biggr )
 =\sum_{k=1}^n v_kx_k^*.
\end{align*}
\end{pf}
\qed

Let $(F, \|\cdot\|_F)$ be a normed Riesz space of real sequences
(cf. [8, p.6] for definition). Following [6], we say that $F$ has
the property $(PS)$, if for all $x\in F$, $x^*$ exists and the
following property holds:
$$
y_1^*+\cdots+y_n^*\le x_1^*+\cdots+x_n^*\quad(n\ge 1)\quad
\Longrightarrow \quad y\in F\mbox{ and } \|y\|_F\le \|x\|_F.
\leqno{(3.2)}
$$
Clearly, for $x\in F$, we have $\tilde x\in F$ and $\|\tilde
x\|_F=\|x\|_F$, where $\tilde x$ is any sequence with $\tilde
x^*=x^*$. In particular, $\tilde x$ can be $x^*$ or any sequence
obtained from $x$ by reordering $x_k$. We have $x_1+\cdots+x_n\le
x_1^*+\cdots+x_n^*$, so the condition in $(3.2)$ can be replaced
by $y_1^*+\cdots+y_n^*\le x_1+\cdots+x_n$. Applying Lemma 3.1, we
get the first main result as follows.

\begin{thm} Let  $(E, \|\cdot\|_E)$, $(F, \|\cdot\|_F)$ be two
normed Riesz space of real sequences with property $(PS)$. In
addition, the following property is satisfied:
$$
  \|x\|_F=\lim_{n\to\infty} \|P_nx\|_F\qquad(x\in F),\leqno{(3.3)}
  $$
  where $P_nx$ is the projection of $x$ onto the first $n$ terms.
Let $A=(a_{j,k})_{j,k\ge 1}$ define an operator from $E$ to $F$,
given by $Ax=y$, where $a_{j,k}\ge 0$ for all $j$ and $k$. If
$(1.3)$ is true for some positive integer $n_0$, then
$\|Ax^*\|_F\ge \|Ax\|_F$ for all $x\in E$ with $x\ge 0$. Hence,
decreasing elements $x$ in $E$ are sufficient to determine
$\|A\|_{E,F}.$
\end{thm}

\begin{pf}
Let $x\in E$ with $x\ge 0$. Then the $(PS)$ property of $E$
implies $x^*\in E$. Since $A$ sends $E$ to $F$, we know that $Ax,
Ax^*\in F$. We claim that $\|Ax^*\|_F\ge \|Ax\|_F$. Let $n\ge
n_0$. By Lemma 2.3, we can find $B=(b_{j,k})\in {\cal
R}_{A_n}^{\gamma,\lambda}$ with $0\le\gamma\le\lambda\le n$ such
that $\{\sum_{k=1}^\infty b_{jk}x_k\}_{j=1}^n$ is decreasing. Let
$C=(c_{j,k})\in {\cal R}_{A_n}$ be the corresponding matrix
obeying $(1.3)$. Let $l$ be fixed. Since $
 \sum_{k=1}^r\left(\sum_{j=1}^lc_{jk}\right)\ge
 \sum_{k\in N_r}\left(\sum_{j=1}^lb_{jk}\right)$
for all $r\ge 1$ and for all $N_r$, it follows from Lemma 3.1 that
$$
 \sum_{k=1}^m \left(\sum_{j=1}^lc_{jk}\right)x_k^*\ge
\sum_{k=1}^m\left(\sum_{j=1}^lb_{jk}\right)x_k\qquad (m\ge 1).
$$
Let $m\to\infty$ and reorder the above sums. Then we obtain
$$
\sum_{j=1}^l \left(\sum_{k=1}^\infty c_{jk}x_k^*\right)
\ge\sum_{j=1}^l\left(\sum_{k=1}^\infty
b_{jk}x_k\right)\qquad(l=1,\ldots,n). \leqno{(3.4)}
$$
For $1\le j\le n$, set $y_j=\sum_{k=1}^\infty c_{jk}x_k^*$ and
$z_j=\sum_{k=1}^\infty b_{jk}x_k$. We also set $y_j=z_j=0$ for
$j>n$. Since $\{z_j\}_{j=1}^\infty$ is decreasing, $z_j^*=z_j$ for
all $j$ and consequently, $(3.4)$ can be rewritten in the form:
$\sum_{j=1}^l y_j\ge\sum_{j=1}^l z_j^*$ for all $l\ge 1$. On the
other hand, $P_n Ax^*\in F$ and it is of the form: $P_n
Ax^*=\{y_1',\cdots,y_n',0,\cdots\}$. Since $C\in {\cal R}_{A_n}$,
$y=\{y_1,\cdots,y_n,\cdots\}$ can be obtained from the sequence
$\{y_1',\cdots,y_n',0,\cdots\}$ by reordering the first $n$ terms.
The $(PS)$ property of $F$ implies $y\in F$, and therefore,
$z=\{z_1,z_2,\cdots\}\in F$. Moreover, $\|P_n Ax^*\|_F=\|y\|_F\ge
\|z\|_F$. We have $B\in {\cal R}_{A_n}$. The same argument as
above also ensures that $\|z\|_F=\|P_nAx\|_F$. Hence,
$\|P_nAx^*\|_F\ge \|P_nAx\|_F$. Let $n\to\infty$. Then $(3.3)$
implies $\|Ax^*\|_F\ge \|Ax\|_F$. Next, consider the case $x\in
E$. Set $\tilde x=\{\tilde x_k\}_{k=1}^\infty$, where $\tilde
x_k=|x_k|$. Then $\tilde x\in E$. Moreover, $\tilde x\ge 0$ and
$\tilde x^*=x^*\in E$. By the result we have proved,
$\|Ax^*\|_F=\|A\tilde x^*\|_F\ge \|A\tilde x\|_F$. We have $
|\sum_{k=1}^\infty a_{j,k}x_k|\le \sum_{k=1}^\infty a_{j,k}|x_k|$
for all $j$. Since $F$ is a normed Riesz space, $\|Ax\|_F\le
\|A\tilde x\|_F$, and consequently, $\|Ax^*\|_F\ge \|Ax\|_F$. This
ensures the validity of $\|A\|_{E,F}=\|A\|_{E,F,\downarrow}.$ We
complete the proof. \qed
\end{pf}

Theorem 3.2 generalizes [6, Theorem 2] and [5, Lemma 2.4]. We
shall investigate them in \S4.

Following the above proof, we see that Theorem 3.2 still holds for
the case of complex sequences, if in addition, elements in $E\cup
F$ satisfy  $\|\tilde x\|=\|x\|$, where
 $\tilde x=\{|x_1|, |
x_2|,\cdots\}$ and $\|\cdot\|$ denotes the corresponding norm in
$E$ or in $F$. Moreover, the assumption that $A$ sends $E$ to $F$
can be removed from Theorem 3.2, whenever $\|Ax^*\|_F$ and
$\|Ax\|_F$ make sense and satisfy
$$\|Ax^*\|_F=\lim_{n\to\infty}\|P_nAx^*\|_F\quad\mbox{ and}\quad
\|Ax\|_F=\lim_{n\to\infty} \|P_nAx\|_F.$$
 In particular, the spaces $E$ and $F$ in Theorem 3.2 can
be one of $\ell_p\quad(1\le p<\infty)$ or $d(w,p)\quad(1\le
p\le\infty)$. However, Theorem 3.2 can not apply to the case
$F=\ell_\infty$ (or $\ell_\infty(w)$), in general. A
counterexample is given by $a_{2,2}=1$, $a_{j,k}=0$ otherwise,
$x_2=1$, and $x_k=0$ for $k\neq 2$. For this example,
$\|Ax^*\|_\infty<\|Ax\|_\infty$ and
$\|A\|_{\ell_p,\ell_\infty}\neq
\|A\|_{\ell_p,\ell_\infty.\downarrow}$ for $1\le p<\infty$. In the
following, we show that $F$ can be $\ell_p(w)$, where $1\le
p<\infty$. Since the case $F=\ell_2(w)$ with $w_n=1/n^3$ fails the
property $(PS)$, Theorem 3.3 is not a special case of Theorem 3.2.

\begin{thm} Let $1\le p<\infty$, $A=(a_{j,k})_{j,k\ge 1}\ge 0$,
and $(E, \|\cdot\|_E)$ be a normed Riesz space of real sequences
with property $(PS)$. If $(1.3)$ is true for $n_0=1$, then
$\|Ax^*\|_{\ell_p(w)}\ge \|Ax\|_{\ell_p(w)}$ for all $x\in E$ with
$x\ge 0$ . Hence, decreasing, non-negative elements $x$ in $E$ are
sufficient to determine  $ \|A\|_{E,\ell_p(w)}.$
\end{thm}

\begin{pf}
We only need the following simple remark and for the same proof in
theorem 3.2 apply: if $x,y \geq 0$ and $\sum_{j=1}^nx_j^p \leq
\sum_{j=1}^ny_j^p$ for all n, then
$\|x\|_{\ell_p(w)}\leq\|y\|_{\ell_p(w)}$. (By Abel summation)
\end{pf}

\section{Investigation of $(1.3)$}

In Theorems 3.2-3.3, we point out that $(1.3)$ is a sufficient
condition for $A$ to guarantee the validity of the equality
$\|A\|_{E, F}=\|A\|_{E, F,\downarrow}$. The purpose of this
section is to find those conditions which are stronger than
$(1.3)$. First, we deal with conditions of Jameson type, that is,
$(1.2)$ and its related conditions. Set $c_{j,k}=a_{j,k}$. We see
that $(1.2)\Longrightarrow (1.3)$. Here $(1.3)$ is assumed for
$n_0=1$. Moreover, the entries of $A^t$ still satisfy $(1.2)$, if
the entries of $A$ do. Here $A^t$ is the transpose of $A$. Hence,
Theorems 3.2-3.3 have the following consequence.

\begin{thm} Theorems 3.2-3.3 remain true, if $(1.3)$ is replaced by $(1.2)$.
Moreover, the conclusions of Theorems 3.2-3.3 also hold for $ A^t$
in place of $A$.
\end{thm}

Clearly, Theorem 4.1 extends [6, Theorem 2] from $E=F$ to any pair
$(E, F)$. Moreover, it can apply to the case $F=\ell_p(w)$, (see
Theorem 3.3), but [6, Theorem 2] fails to do so. We know that
$(4.1)\Longrightarrow (1.2)$:
$$
 a_{j,k}\ge a_{j+1,k}\quad(j,k\ge 1)\quad\mbox{ and}\quad
    \sum_{j=1}^la_{j,k}\ge\sum_{j=1}^la_{j,k+1}\quad(k,l\ge 1),\leqno{(4.1)}
    $$
(see [6, Proposition 3]). Hence, Theorem 4.1 has  the following
consequence.

\begin{cor} Theorems 3.2-3.3 remain true, if $(1.3)$ is replaced by $(4.1)$.
Moreover, the conclusions of Theorems 3.2-3.3 also hold for $ A^t$
in place of $A$.
\end{cor}

In [6, Proposition 3], Jameson pointed out that the following
condition also implies $(1.2)$, and so Corollary 4.2 still holds,
if we replace $(4.1)$ by $(4.1^*)$:
$$
 a_{j,k}\ge a_{j,k+1}\quad(j,k\ge 1)\quad\mbox{ and}\quad
    \sum_{k=1}^ra_{j,k}\ge\sum_{k=1}^ra_{j+1,k}\quad(j,r\ge 1).\leqno{(4.1^*)}
    $$
    We shall prove in Corollary 4.7 that the second condition in
    $(4.1^*)$ is redundant.

 Corollary 4.2 can apply to the Hilbert matrix $H=(h_{j,k})_{j,k\ge 1}$, defined by $h_{j,k}=1/(j+k-1)$. It can also
 apply to the weighted mean matrix
$A_W^{WM}=(a_{j,k}^{WM})_{j,k\ge 1}$ and the N\"orlund mean matrix
$A_W^{NM}=(a_{j,k}^{NM})_{j,k\ge 1}$, where
$a_{j,k}^{WM}=a_{j,k}^{NM}=0$ for $j<k$ and
$$
 a_{j,k}^{WM}=w_k/(w_1+\cdots+w_j)\quad(j\ge k),\leqno{(4.2)} $$
$$
 a_{j,k}^{NM}=w_{j-k+1}/(w_1+\cdots+w_j)\quad(j\ge k).\leqno{(4.3)}
$$

\begin{cor} Let $w_1>0$ and $w_n\ge 0$ for all $n>1$. Then Theorems 3.2-3.3 remain true,
if $(1.3)$ is replaced by any of $(i)$ and $(ii)$:
\begin{itemize}
\item[$(i)$] $A=(A_W^{WM})^t$ with $w_n\downarrow$.
\item[$(ii)$] $A=(A_W^{NM})^t$, where $w_n\uparrow$ and $w_{n+1}/w_n\le w_n/w_{n-1}$ for all $n$.
\end{itemize}
Moreover, the conclusions of Theorems 3.2-3.3 also hold for $ A^t$
in place of $A$.
\end{cor}

\begin{pf}
Obviously, $(A_W^{WM})^t\ge 0$ and $(A_W^{NM})^t\ge 0$. Consider
Case $(i)$. Set $A=(A_W^{WM})^t=(a_{j,k})_{j,k\ge 1}$. It is easy
to see that $a_{j,k}\ge a_{j+1,k}$ for all $j,k\ge 1$ if and only
if $w_n\downarrow$. Moreover, we have
$$
  \sum_{j=1}^l a_{j,k}=\left \{
 \begin{array}{ll}
  \frac {w_1+\cdots+w_l}{w_1+\cdots+w_k}&\qquad (l\le k),\\
  1     &\qquad (l>k).
\end{array} \right.
$$ This implies $\sum_{j=1}^la_{j,k}\ge\sum_{j=1}^la_{j,k+1}$ for
all $k,l\ge 1$. The above argument shows that $(4.1)$ holds. By
Corollary 4.2, we get $(i)$. Next, consider $(ii)$. It is an
consequence of the following well-known lemma. Write $A_n=
\sum_{j=1}^n a_j$ and $B_n$ similarly. If $(a_n/b_n)$ is
increasing(or decreasing), then so is $A_n/B_n$. This shows that
$
 \sum_{j=1}^l a_{j,k+1}\le \sum_{j=1}^l a_{j,k}$
 for $l<k$.
This inequality is also true for the case $l\ge k$, because
$\sum_{j=1}^l a_{j,k+1}\le 1=\sum_{j=1}^l a_{j,k}.$ Thus, $(4.1)$
is satisfied. By Corollary 4.2, we get $(ii)$. This completes the
proof.\qed
\end{pf}

The conclusion of Corollary 4.3(i) for $A^t$ and  for $E=F=\ell_p$
was established by Bennett in [2, page 422] and [3, page 422],
where $1<p<\infty$. For $w_n=\binom {n+\alpha-2}{n-1}$, $A_W^{WM}$
and  $A_W^{NM}$ are denoted by $\Gamma(\alpha)$ and $C(\alpha)$,
respectively. They are called  the Gamma matrix and the Ces\`aro
matrix, of order $\alpha$, (cf. [3, p.410], [4] \& [9, Chapter
III]). We know that $w_n\uparrow$ for $\alpha\ge 1$ and
$w_n\downarrow$ for $0<\alpha\le 1$ (cf. [9, page 77]). Moreover,
for $\alpha\ge 1$, we have
$$
 \frac {w_{n+1}}{w_n}=\frac {n+\alpha-1}n\le
 \frac {n+\alpha-2}{n-1}= \frac {w_n}{w_{n-1}}.
$$
Hence, by Corollary 4.3, the conclusions of Theorems 3.2-3.3 hold
for $A$ to be any of the matrices: $\Gamma(\alpha),
\Gamma(\alpha)^t\quad(0<\alpha\le 1)$ and $C(\alpha),
C(\alpha)^t\quad(\alpha\ge 1)$.

Following [3], we say that $A=(a_{j,k})_{j,k\ge 1}$ is a
summability matrix, if $A$ is a non-negative lower triangular
matrix with $\sum_{k=1}^\infty a_{j,k}=1$ for all $j$. For such
type of matrices, we have the following result.

\begin{cor} Let $A=(a_{j,k})_{j,k\ge 1}$ be a summability matrix.
Then $(4.4)\Longrightarrow (1.2)$, where
$$
a_{j,k}\ge \max (a_{j+1,k}, a_{j+1,k+1})\qquad(j\ge k\ge
1).\leqno{(4.4)}
$$
Hence, Theorems 3.2-3.3 remain true, if $(1.3)$ is replaced by
$(4.4)$. Moreover, the conclusions of Theorems 3.2-3.3 also hold
for $ A^t$ in place of $A$.
\end{cor}

\begin{pf}
The second part follows from Theorem 4.1. We claim that
$(4.4)\Longrightarrow (1.2)$. Divide the proof into three cases.
Case I is $l\le r$. For this case, we have
$$
 \sum_{j=1}^l\sum_{k=1}^r a_{j,k}\ge \sum_{j=1}^l\sum_{k=1}^l a_{j,k}=l.\leqno{(4.5)}
$$
On the other hand, we know that $A$ is a summability matrix. Thus,
$ \sum_{k=1}^\infty a_{j,k}=1$ for all $j$. This implies
$$
 \sum_{j\in N_l}\sum_{k\in N_r} a_{j,k}\le
 \sum_{j\in N_l}\biggl (\sum_{k=1}^\infty a_{j,k}\biggr )=|N_l|=l.\leqno{(4.6)}
$$
Putting $(4.5)-(4.6)$ together yields $(1.2)$ for Case I. Next,
consider the case: $l>r$ and $N_r=\{1,2,\cdots,r\}.$ Write
$N_l=\{j_1,\cdots,j_l\}$ in the alphabet order. Then
$$
\sum_{s=1}^r\sum_{k\in N_r} a_{j_s,k} \le \sum_{s=1}^r\biggl
(\sum_{k=1}^\infty a_{j_s,k}\biggr ) =r=\sum_{s=1}^r\sum_{k=1}^r
a_{s,k}.\leqno{(4.7)}
$$
On the other hand, for $r<s\le l$ and $k\in N_r$, by $(4.4)$, we
get $a_{s,k}\ge a_{j_s,k}$, and so $ \sum_{k\in N_r} a_{j_s,k}\le
\sum_{k=1}^r a_{s,k}.$ Sum up both sides over $s\in \{r+1,\cdots,
l\}$. Then
$$
 \sum_{s=r+1}^l\sum_{k\in N_r} a_{j_s,k}
\le \sum_{s=r+1}^l\sum_{k=1}^r a_{s,k}.\leqno{(4.8)}$$ Putting
$(4.7)-(4.8)$ together yields $
 \sum_{j\in N_l}\sum_{k\in N_r} a_{j,k}
\le \sum_{s=1}^l\sum_{k=1}^r a_{s,k}. $ This is $(1.2)$. It
remains to prove the case that $l>r$ and $N_r$ is any set of
positive integers with $|N_r| =r$. Write $N_r=\{k_1,\cdots,k_r\}$
in the alphabet order. We can assume that $j_1\ge k_1$, otherwise,
$ a_{j_1,k}=0$ for all $k\in N_r$. In this case,  $ \sum_{k\in
N_r}a_{j_1,k}=0$, which allows us to replace $a_{j_1,k}$, with
$k\in N_r$, by $a_{j,k}$ for some $j\notin N_l$. Let $N_l$ be the
corresponding new index set. Our replacement leads us to deal with
a case, which has a bigger sum on the right side of $(1.2)$.
Similarly, we can assume $j_l\ge k_r$. With the help of $(4.4)$,
we can replace $a_{j_s,k_t}$ by $a_{j_s-k_1+1,k_t-k_1+1}$. After
this replacement, we can assume $k_1=1$. We shall prove that under
suitable replacements, we can assume $k_t=t$ for all $t=2,\cdots,
r$. For any $t^*$ with $k_{t^*+1}\ge k_{t^*}+2$, let $s^*$ be the
smallest integer with $j_{s^*}\ge k_{t^*+1}$. This $s^*$ exists,
because $j_l\ge k_r\ge k_{t^*+1}$. If $s^*>1$, then
$j_{s^*-1}<k_{t^*+1}$ and so $a_{j_s,k_t}=0$ for all $1\le s<s^*$
and $t^*<t\le r$. Here we use the fact that $A$ is a lower
triangular matrix. Replace $a_{j_s,k_t}$ by $a_{j_s-1,k_t-1}$ for
$s^*\le s\le l$ and $t^*<t\le r$. Simultaneously, we make the
change $a_{j_s,k_t}\longrightarrow a_{j_s-1,k_t}$ for $s^*\le s\le
l$ and $1\le t\le t^*$, whenever $j_s>j_{s-1}+1$. If $s^*=1$, then
$j_{s^*}\ge k_{t^*+1}\ge 2$. Replace $a_{j_s,k_t}$ by
$a_{j_s-1,k_t-1}$ with $s^*\le s\le l, t^*<t\le r$, and make the
change $a_{j_s,k_t}\longrightarrow a_{j_s-1,k_t}$ for all $s^*\le
s\le l$ and $1\le t\le t^*$. The above argument shows that the
difference $k_{t^*+1}-k_{t^*}$ can be reduced by 1 for each
replacement. Continue this process several times and we finally
get $k_{t^*+1}=k_{t^*}+1$. Our argument leads us to the choice
$k_t=t$ for all $t$ and our problem reduces to Case II. This has
been proved before. Hence, the desired result follows.\qed
\end{pf}

Corollary 4.4 allows us to deal with the case $A=A_W^{NM}$ with
$w_n\downarrow$.

\begin{cor} Let $w_1>0$ and $w_n\ge 0$ for all $n>1$. Then Theorems 3.2-3.3 remain true,
if $(1.3)$ is replaced by $A=A_W^{NM}$ with $w_n\downarrow$.
Moreover, the conclusions of Theorems 3.2-3.3 also hold for $
(A_W^{NM})^t$ in place of $A$.
\end{cor}

\begin{pf}
We know that $A_W^{NM}=(a_{j,k}^{NW})_{j,k\ge 1}$ is a summability
matrix. The hypothesis that $w_n\ge 0$ and $w_n\downarrow$ implies
$$
 \frac {w_1+\cdots+w_{j+1}}{w_1+\cdots+w_j}
 \ge 1\ge \frac {w_{j-k+2}}{w_{j-k+1}}\qquad (j\ge k\ge 1).$$
This leads us to $(4.4)$ for $a_{j,k}=a_{j,k}^{NW}$. By Corollary
4.4, we get the desired result.\qed
\end{pf}

 For $w_n=\binom {n+\alpha-2}{n-1}$, $A_W^{NM}=C(\alpha)$. Moreover,
 $w_n\downarrow\Longleftrightarrow 0<\alpha\le 1$. Hence, by Corollary
 4.5, the conclusions of Theorems 3.2-3.3 hold for $A$ to be one of the matrices:
 $C(\alpha), C(\alpha)^t\quad(0<\alpha\le 1)$.

The matrix $A_W^{NM}$ involved in Corollary 4.5 is row increasing
in the triangular sense, that is, $a_{j,k}\le a_{j,k+1}$ for all
$j>k$. This fact does not imply that Corollary 4.5 can be extended
to any summability matrix with rows increasing in the triangular
sense. A counterexample is given below: $$
 A=\begin{pmatrix}
  1      & 0      & 0      & 0      & 0       & \cdots \\
  0      & 1      & 0      & 0      & 0       & \cdots \\
  0      & 1/2    & 1/2    & 0      & 0       & \cdots \\
  0      & 0      & 0      & 1      & 0       & \cdots \\
  0      & 0      & 0      & 0      & 1       & \cdots \\
  \vdots & \vdots & \vdots & \vdots & \vdots  & \ddots
\end{pmatrix}\;.
$$ For $p=1$, we have $\|A\|=\|Ae_2\|_1=3/2$, but for decreasing
$x=(x_n)$,
\begin{align*}
 \mbox{(4.9)\hskip 0.5in}\|Ax\|_1
 &= x_1+\frac{3}{2}x_2+\frac{1}{2}x_3+x_4+\cdots\\
  &\leq \frac{5}{4}(x_1+x_2)+\cdots \leq \frac{5}{4}\|x\|_1.
\end{align*}

In Theorem 4.1, we deal with the condition $(1.2)$, which
corresponds to the case $c_{j,k}=a_{j,k}$ of $(1.3)$. In the
following, we consider the case $c_{j,k}=b_{j,k}$. More precisely,
we consider the following condition for $n\ge n_0$:
$$
 \sum_{j=1}^l\sum_{k=1}^r b_{j,k}\ge
 \sum_{j=1}^l\sum_{k\in N_r} b_{j,k}
 \qquad(1\le l\le n;\quad r\ge 1;\quad |N_r|=r),\leqno{(4.10)}
$$
where $n_0$ and $B=(b_{j,k})$ are stated in $(1.3)$. We know that
$(4.10)\Longleftrightarrow (4.10^*)$:
$$
 \biggl \{\sum_{j=1}^l b_{j,k}\biggr \}_{k=1}^\infty
 \quad\mbox{is a decreasing sequence}
   \qquad(1\le l\le n).\leqno{(4.10^*)}
$$
This fact can be derived by considering $N_r=\{1,\cdots, r-1,
r+1\}$. By Theorems 3.2-3.3, we obtain the following result.

\begin{thm} Theorems 3.2-3.3 remain true, if $(1.3)$ is replaced by any of $(4.10)$ and $(4.10^*)$.
\end{thm}

The matrix $A=(a_{j,k})_{j,k\ge 1}$, with
$a_{1,1}=a_{2,2}=a_{2,3}=1$ and 0 otherwise, obeys the inequality:
$
 \|A\|_{\ell_2,\ell_2}>\|A\|_{\ell_2,\ell_2,\downarrow}$. This follows from
 the fact that
$$
 \|A\|_{\ell_2,\ell_2}=\sup_{\|x\|_2=1, x\ge 0} (x_1^2+x_2^2+x_3^2+2x_2x_3)^{1/2}$$
 is attained only at
 $\displaystyle
 x=(0,\frac 1{\sqrt 2}, \frac 1{\sqrt 2}, 0,\cdots)$, which is not a decreasing sequence.
 This example shows that Theorem 4.6 is not true, if $(4.10)$ is replaced by $(4.11)$:
$$
 \sum_{j=1}^l\sum_{k=1}^r a_{j,k}\ge
 \sum_{j=1}^l\sum_{k\in N_r} a_{j,k}
 \qquad(l, r\ge 1;\quad |N_r|=r).\leqno{(4.11)}
$$
Obviously, $(4.11)$ is weaker than $(1.2)$. We know that $(4.11)$
is equivalent to the second part of $(4.1)$. Hence, the first part
of $(4.1)$ can not be removed from Corollary 4.2.

 For $\displaystyle
B=(b_{j,k})\in\cup_{0\le\gamma\le\lambda\le n}
 {\cal R}_{A_n}^{\gamma,\lambda}$,
$\sum_{j=1}^l b_{j,k}=\sum_{j\in N_l} a_{j,k}$ for some index set
$N_l$ with $|N_l|=l$. Hence, $(4.12)\Longrightarrow (4.10^*)$:
\begin{itemize}
\item[(4.12)] $\hskip 0.2in A$ is row decreasing, that is,
    $a_{j,k}\ge a_{j,k+1}$ for all $j,k\ge 1$.
\end{itemize}
As a consequence of Theorem 4.6, we obtain the following result.

\begin{cor} Theorems 3.2-3.3 remain true, if $(1.3)$ is replaced by
$(4.12)$.
\end{cor}

Corollary 4.7 extends [5, Lemma 2.4] from the pair
$(\ell_p,\ell_p)$ to the pair $(E, F)$. Moreover, it indicates
that the condition $(5^*)$ in [6, Proposition 3] is enough to
ensure the validity of [6, Theorem 2]. Obviously, the entries of
the Hilbert matrix $H$ satisfy $(4.12)$. Hence, the conclusions of
Theorem 3.2-3.3 hold for $A=H$. Applying Corollary 4.7 to the
N\"orlund mean matrix $A_W^{NW}$, we get the following
consequence.

\begin{cor} Let $w_1>0$ and $w_n\ge 0$ for all $n>1$. Then Theorems 3.2-3.3 remain true,
if $(1.3)$ is replaced by $A=A_W^{NM}$ with $w_n\uparrow$.
\end{cor}

Corollary 4.8 is a generalization of Corollary 4.3(ii) for the
N\"orlund mean matrix $A_W^{NM}$. For this matrix, the condition
$w_{n+1}/w_n\le w_n/w_{n-1}$, required in Corollary 4.3(ii), is
redundant. However, we do not know whether this condition can be
removed for the transpose $(A_W^{NM})^t$. For the case $w_n=\binom
{n+\alpha-2}{n-1}$, it does, (see the statement given after the
proof of Corollary 4.3). It is still open for general $w_n$.

\end{document}